\documentclass[11pt]{article}
\usepackage[latin1]{inputenc}
\usepackage[T1]{fontenc}
\usepackage{amsmath}
\usepackage{bbm}
\usepackage{amsthm}
\usepackage{amssymb,mathrsfs}
\usepackage{fourier}
\usepackage[english]{babel} 
\usepackage[all]{xy}
\usepackage{setspace}
\usepackage{color}
\usepackage{tabularx}
\usepackage{amsmath}
\newtheorem{theorem}{Theorem}[section]
\newtheorem{corollary}[theorem]{Corollary}

\newtheorem{lemma}[theorem]{Lemma}
\newtheorem{proposition}[theorem]{Proposition}

\newcommand{\sinc}{sinc}
\allowdisplaybreaks[1]

\newcommand\N{{\mathbb N}}
\newcommand\Z{{\mathbb Z}}

\newcommand\R{{\mathbb R}}
\newcommand\C{{\mathbb C}}
\newcommand\T{{\mathbb T}}

\begin{document}

\title{Bounded cosine functions close to continuous scalar bounded cosine functions}

\date{\relax}
\author{Jean Esterle}

\maketitle

{\bf {Abstract}}: Let $(C(t))_{t\in \R}$ be a cosine function in a unital Banach algebra. We show that if $sup_{t\in \R}\Vert C(t)-c(t)\Vert <2$ for some continuous scalar bounded 
cosine function $(c(t))_{t\in \R},$ then the closed subalgebra generated by $(C(t))_{t\in \R}$ is isomorphic to $\C^k$ for some positive integer $k.$ If, further, $sup_{t\in \R}\Vert C(t)-c(t)\Vert <{8\over 3\sqrt 3},$ then $C(t)=c(t)$ for $t\in \R.$ 

{\it Keywords: Cosine function, scalar cosine function, commutative local Banach algebra.

}
\section{Introduction}

 Recall that a cosine function taking values in  a unital normed algebra $A$ with unit element $1_A$ is a family $C=(C(t))_{t\in \R}$ of elements of $A$  satisfying the so-called d'Alembert equation

\begin{equation}C(0)=1_A, C(s+t)+C(s-t)=2C(s)C(t)  \ \ (s\in \R,  t\in \R).\end{equation}

Equality (1) is also used to define $G$-cosine families $C=(C(g))_{g\in G}$ over an abelian group $G,$ and in particular cosines sequences $(C_n)_{n\in \Z}.$

A cosine function $C= (C(t))_{t\in \R}$ is said to be bounded if there exists $M>0$ such that $\Vert C(t)\Vert \le M$ for every $t \in \R.$ In this case we set

$$\Vert C \Vert_{\infty} =sup_{t \in R}\Vert C(t)\Vert, \ dist(C_1,C_2)=\Vert C_1-C_2\Vert_{\infty}.$$
A  cosine function is said to be scalar if $C(t)\in \C.1_A$ for $t \in \R.$   It is well-known and easy to see that a bounded  complex-valued cosine function $c$ takes values in $[-1,1]$ and that a bounded continuous complex-valued cosine function satisfies $c(t)=cos(at)$ where 
$a=2lim_{t\to 0}{1-c(t)\over t^2} \in \R.$  More generally if $X$ is a Banach space, and if $C=(C(t))_{t\in \R}$  is a strongly continuous ${\mathcal B}(X)$-valued cosine function, then the generator $a$ of $C$ is bounded if and only if the cosine function $C$ is continuous at $0$ with respect to the operator norm on ${{\mathcal B}(X)},$ and in this situation $C(t)=cos(ta),$ where $cos(ta)$ is defined by the usual series, see for example \cite{so}.

Strongly continuous operator valued cosine functions are a classical tool in the study of differential equations, see for example \cite{abhn}, \cite{beh}, \cite{n}, \cite{tw}, and a functional calculus approach to these objects was  developped recently in \cite{h1}. 

Bobrowski and Chojnacki proved recently in \cite{bc}  that if the set of all bounded strongly continuous cosine
functions on a Banach space $X$ is treated as a metric space under the operator norm on ${\mathcal B}(X),$
then the isolated points of this set are precisely the scalar cosine
functions.

They also gave a quantitative version of this result, by showing that if a strongly continuous operator valued cosine function on a Banach space $C(t)$ satisfies $sup_{t\ge 0}\Vert C(t) -c(t)\Vert <1/2$ for some scalar bounded continuous cosine function $c(t)$ then $C(t)=c(t)$ for $t \in \R,$ and  Schwenninger and Zwart showed in \cite{zs} that this result remains valid under the condition $sup_{t\ge 0}\Vert C(t) -c(t)\Vert <1.$ An elementary proof of this result is given by Chojnacki in \cite{c}.

 The purpose of this paper is to show that this result holds when  $sup_{t\ge 0}\Vert C(t) -c(t)\Vert <{8\over 3\sqrt 3},$  which is the optimal constant since $sup_{t\in \R}\vert cos(at)-cos(3at)\vert ={8\over 3\sqrt 3}$ for every $a\neq0,$ and that no continuity condition on $C$ is needed if the scalar bounded cosine function $c$ is assumed to be continuous, see theorem 3.6 (ii). In fact, this $"0-{8\over 3\sqrt 3}$ law" was already proved in a very recent paper by Bobrowski, Chojnacki and Gregoriewicz \cite{bcg}, which appeared after the present paper was submitted. The methods developed here vary in various aspects from those used by
Bobrowski et al.
 
 We also point out that the 'generic' distance between two continuous real-valued bounded cosine functions is 2, a consequence of the fact that finite independent subsets of the torus are "Kronecker sets": if $a\neq 0,$ then the set $\Omega(a,m):=\{ b\ge 0 \ \ sup_{t\in \R} \vert cos(at)-cos(bt)\vert \le m\}$ is finite, and every $b\in \Omega(a,m)$ has the form $b={pa\over q},$ where $p,q$ are odd, $gcd(p,q)=1$, $1\le p \le {\pi\over arccos(m-1)}, 1\le q \le {\pi\over arccos(m-1)}.$ 
 
 This description of the set $\Omega(a,m)$ leads to a description of cosine functions $C=(C(t))_{t\in \R}$ in a Banach algebra $A$ satisfying $sup_{t\in \R}\Vert C(t) -c(t)\Vert =m$ when $m<2.$ In this case we show in theorem 3.6 (i) that there exists $k \le card \left ( \Omega(a,m)\right )$ such that the closed subalgebra $A_1$ generated by $C$ is isomorphic to $\C^k,$ and we also show that there exists a family $p_1,\dots, p_k$ of pairwise orthogonal idempotents of $A_1$ and a family $(b_1,\dots,b_k)$ of distinct elements of $\Omega(a,m)$ such that we have
 
 $$C(t)=\sum \limits_{j=1}^kcos(b_jt)p_j.$$
 
 This implies in particular that if a cosine family $(C(t))_{t\in \R}$ satisfies $sup_{t\in \R}\Vert C(t)-1_A\Vert <2,$ then $C(t)=1_A$ for $t \in \R.$ This result was proved very recently
 by Schwenninger and Zwart in \cite{zs1} for strongly continuous cosine families of bounded operators, but the general case seems new. 
 
The description of the set $\Omega(a,m)$ pertains to folklore, but the operator theoretical part of the proofs seems new. It is based on the fact that every bounded  cosine function $(C(t))_{t\in \R}$ taking values in a commutative unital Banach algebra having a unique maximal ideal is scalar, see
theorem 2.3 and corollary 2.4 in section 2.

Since $sup_{n\ge 1}\left \vert 1 -cos\left ({2n\pi\over 3}\right )\right \vert={3\over 2},$ the constant ${8\over 3\sqrt 3}$ does not work for cosine sequences. It is nevertheless possible to show that if a $G$-cosine family $C$ in a unital Banach algebra $A$ satisfies $sup_{g\in G}\left \Vert C(g)-c(g)\right \Vert <{\sqrt 5\over 2}$ for some bounded scalar $G$-cosine family $c,$ then $C(g)=c(g)$ for $g\in G.$ and the constant ${\sqrt 5\over 2}= max_{n\ge 1}\left | cos \left ( {n\pi \over 5} \right )-cos\left ( {3n\pi \over 5} \right )\right |$ is obviously optimal. Details will be given elsewhere.

The author would like to give his very  warm thanks to the referee for his insightful comments and corrections.


\section{Cosine sequences in commutative  local Banach algebras}

Set $f(x)=arccos(x).$ Then $f'(x)=-\frac{1}{\sqrt{1-x^2}},$ and we have, for $x\in ]-1,1[$

$$\frac{1}{\sqrt{1-x^2}}=1+\sum\limits_{n=1}^{+\infty}\frac{(-1)^n}{n!}\left(-\frac{1}{2}\right )\left (- \frac{1}{2}-1\right )\dots \left (- \frac{1}{2}-n+1\right )x^{2n}$$ $$=
1+\sum \limits_{n=1}^{+\infty}\frac{(2n)!}{2^{2n}n!^2}x^{2n}.$$

Hence we have, for $x\in (-1,1),$ with the convention $0!=1,$

$$arccos(x)=\frac{\pi}{2} -\sum\limits_{n=0}^{+\infty}\frac{(2n)!}{2^{2n}(2n+1)n!^2}x^{2n+1}.$$

It follows from example from Stirling's formula and Riemann's criterion that the series $\sum\limits_{n=0}^{+\infty}\frac{(2n)!}{2^{2n}(2n+1)n!^2}$ is convergent, and it follows from Abel's lemma that the power series expansion of $arccos(x)$ remains valid for $x=1$ and $x=-1.$ Since $arccos(-1)=\pi,$ we have

$$\sum\limits_{n=0}^{+\infty}\frac{(2n)!}{2^{2n}(2n+1)n!^2}=\frac{\pi}{2}.$$

Now let $A$ be a unital Banach algebra of unit element $1_A.$ We will write $\lambda=\lambda.1_A$ when $\lambda$ is scalar if there is no risk of confusion. We define $e^x$ by the usual series and set, for $x\in A,$

$$cos(x)=\frac{e^{ix}+e^{-ix}}{2}=\sum_{n=0}^{+\infty}(-1)^n\frac{x^{2n}}{(2n)!},$$

so that $cos(x+i\pi n)=(-1)^ncos(x)$ for $n\in \Z.$ When $z\in \C,$ this gives the usual cosine function of a complex variable $z.$

Let $\widehat A$ be the Gelfand space of $A,$ i.e. the space of all algebra homomorphisms from $A$ onto $\C.$ Notice that if $\chi \in \widehat A$ we have

$$\chi(cos(x))=cos(\chi(x)).$$

If sup$_{n\ge 1}\Vert x^n\Vert \le M <+\infty$, then the series $\frac{\pi}{2} -\sum\limits_{n=0}^{+\infty}\frac{(2n)!}{2^{2n}(2n+1)n!^2}x^{2n+1}$ is convergent,
we can set 

$$arccos(x)=\frac{\pi}{2} -\sum\limits_{n=0}^{+\infty}\frac{(2n)!}{2^{2n}(2n+1)n!^2}x^{2n+1},$$

and we have

\begin{equation} \Vert arccos(x)\Vert \le \frac{M+1}{2}\pi, \chi(arccos(x))=arccos(\chi(x)) \ \ (\chi \in \widehat A).\end{equation}

Also it follows from standard properties of the holomorphic functional calculus that $cos(arccos(\lambda x))=\lambda x$ for $\vert \lambda \vert <1.$ By continuity, we obtain the tautological formula

\begin{equation}cos(arccos(x))=x.\end{equation}

\begin{proposition} Let $A$ be a unital Banach algebra, and let $(c_n)_{n\in \Z}\subset A$ be a cosine sequence. If sup$_{n\ge 1}\Vert c_n\Vert \le M <+\infty,$ then sup$_{p\ge 1}\Vert c_n^p\Vert\le M$ for $n\ge 1,$ 
and we have

\begin{equation}\Vert arccos(c_n)\Vert \le \frac{M+1}{2}\pi.\end{equation}

Moreover $Spec(c_n)\subset [-1,1]$ for every $n \in \N,$ and for every character $\chi$ on $A$ we have 

$$\chi(c_n)=cos(n\beta_{\chi}) \ (n\ge 1),$$

where $\beta_{\chi}=\chi(arccos(c_1))=arccos(\chi(c_1))\in [0,\pi].$
\end{proposition}

Proof: Let $p\ge 1,$ and assume that we have

\begin{equation}c_1^p=\sum\limits _{k=0}^p\alpha_{k,p}c_k,\end{equation}

where $\alpha_{k,p} \ge 0,$ $\sum \limits_{k=0}^p\alpha_{k,p}=1,$ which is trivially true for $p=1.$ Using (1), we obtain

$$c_1^{p+1}=\sum_{k=0}^p\alpha_{k,p}c_1c_k = \alpha_{0,p}c_1 +\sum\limits_{k=1}^p{\alpha_{k,p}\over 2}c_{k-1}+\sum\limits _{k=1}^p{\alpha_{k,p}\over 2}c_{k+1}$$ $$\alpha_{0,p}c_1 +\sum\limits_{k=0}^{p-1}{\alpha_{k+1,p}\over 2}c_{k}+\sum\limits _{k=2}^{p+1}{\alpha_{k-1,p}\over 2}c_{k}=\sum \limits_{k=0}^{p+1}\alpha_{k,p+1}c_k,$$

where $$\alpha_{0,p+1}={\alpha_{1,p}\over 2}, \alpha_{1,p+1}=\alpha_{0,p}+{\alpha_{2,p}\over 2}, \alpha_{k,p+1}={\alpha_{k-1,p}+\alpha_{k+1,p}\over 2} \ \mbox{for}\ 2\le k\le p-1,$$
$$\alpha_{p,p+1}={\alpha_{p-1,p}\over 2}, \alpha_{p+1,p+1}={\alpha_{p,p}\over 2}.$$

Clearly, $\alpha_{k,p+1}\ge 0,$ and $\sum\limits _{k=0}^{p+1}\alpha_{k,p}=1.$ We thus see that (5) holds for every $p\ge 1.$
 Hence $sup_{p\ge 1}\Vert c_1^p\Vert \le m.$ Applying this result to the cosine sequence $(c_{nm})_{m\in \Z},$ we see that $sup_{p\ge 1}\Vert c_n^p\Vert \le M$ for every $n\ge 1.$ Inequality (4) follows then from (2).

If $\chi$ is a character on $A,$ then the sequence $(\chi(c_n))_{n\ge 1}$ is a bounded complex-valued cosine sequence. Hence $\chi(c_n) \in [-1,1]$ for $n\ge 1,$ and we have $c_n=cos(n\beta)$ for $n\ge 1,$ where $\beta$ is any real number satisfying $c_1=cos(\beta).$ This holds in particular when $\beta=\beta_{\chi}=\chi(arccos(c_1))=arccos(\chi(c_1))\in [0,\pi].$

 $\square$

Since $cos(arccos(c_1))=c_1,$ we have $cos(narccos(c_1))=c_n,$ and we obtain an alternative approach to the group decomposition  given in \cite{c} by setting $\nu=e^{iarccos(c_1)}.$

Now define the sine function and the "cardinal sine" function on a unital Banach algebra $A$ by the usual formulae

$$sin(x)=\sum_{n=0}^{\infty}(-1)^n\frac{x^{2n+1}}{(2n+1)!}, \sinc(x)=\sum_{n=0}^{\infty}(-1)^n\frac{x^{2n}}{(2n+1)!},$$

so that $\sinc(0)=1,$ $x\sinc(x)=sin(x).$ We have again $\chi(sin(x))=sin(\chi(x))$ and $\chi(\sinc(x))=\sinc(\chi(x))$ for $\chi \in \widehat A.$ For $x,y \in A$ such that $yx=xy,$ we have the usual formula

$$cos(x)-cos(y)=2sin\left ( \frac{y-x}{2}\right )sin\left ( \frac{x+y}{2}\right ).$$

Recall $x \in A$ is said to be quasinilpotent if lim$_{n\to +\infty}\Vert x^n\Vert=0,$ which is equivalent to the fact that $\chi(x)=0$ for every $\chi \in \widehat A$ if $A$ is commutative.

\begin{lemma} Let $A$ be a commutative Banach algebra, let $x\in A$ and $y\in A$ be two quasinilpotent elements of $A$ and let $\lambda \in \C.$ 

(i) If $\lambda \notin \pi\Z,$ and if $cos(\lambda.1_A +x)=cos(\lambda.1_A +y),$ then $x=y.$

(ii) If $\lambda \in \pi \Z,$ and if $cos(\lambda.1_A +x)=cos(\lambda.1_A +y),$ then $x^2=y^2.$

\end{lemma}

Proof: If $cos(\lambda.1_A +x)=cos(\lambda.1_A +y),$ we have

$$sin\left ( \frac{y-x}{2}\right) sin\left (\lambda +\frac{x+y}{2}\right )=0.$$

If $\lambda \notin \pi\Z,$ we have, for $\chi \in \widehat A,$

$$\chi \left (sin\left (\lambda.1_A +\frac{x+y}{2}\right )\right) =sin(\lambda)\neq 0,$$

and so $sin\left (\lambda.1_A +\frac{x+y}{2}\right )$ is invertible and further

$$(x-y)\sinc\left ( \frac{y-x}{2}\right)=2sin\left ( \frac{y-x}{2}\right)=0.$$

Since $\chi \left (\sinc\left ( \frac{y-x}{2}\right)\right )=\sinc(0)=1$ for every $\chi \in \widehat A,$ $\sinc\left ( \frac{y-x}{2}\right)$ is invertible and $x=y.$

If $\lambda \in \pi\Z,$ then we have $cos(x)=cos(y),$ which gives

$$(y^2-x^2)\sinc\left ( \frac{y-x}{2}\right) \sinc\left (\frac{x+y}{2}\right)=4sin\left ( \frac{y-x}{2}\right) sin\left (\frac{x+y}{2}\right )=0.$$

We see again that $\sinc\left ( \frac{y-x}{2}\right)$ and $\sinc\left ( \frac{x+y}{2}\right)$ are invertible, which shows that $x^2=y^2.$

\begin{theorem} Let $(c_n)_{n\in \Z}$ be a bounded cosine sequence in a Banach algebra $A,$ and assume that $spec(c_1)$ is a singleton. Then $(c_n)_{n\ge 1}$ is a scalar sequence, and there exists
$b\in \R$ such that $c_n=cos(nb).1_A$ for $n\ge 1.$
\end{theorem}

Proof: Set $M=sup_{n\ge 1}\Vert c_n\Vert.$ We can assume that $A$ is a commutative Banach algebra generated by $c_1,$ so that $\widehat A$ consists of a single character $\chi.$ Let $\lambda_n=\chi(c_n)$ be the unique element of $spec(c_n),$ and set $\beta=\beta_{\chi}=arccos(\chi(c_1))=\chi(arccos(c_1)),$ so that $\lambda_n=cos(n\beta)$ for $n\ge 1.$

Set $x_n=arccos(c_n).$ Since $c_n=cos(x_n),$ we have
$cos(n\beta)=cos(\chi(x_n)).$ It follows then from standard properties of the cosine function on $\C$ that there exists $k_n\in \Z$ such that $\chi(x_n)=\pm n\beta +2k_n\pi,$
and we have $\chi(x_1)=\beta.$

If $\chi(x_n)=n\beta +2k_n\pi,$ set $y_n=x_n-2k_n\pi.1_A,$ and if $\chi(x_n)=-n\beta +2k_n\pi$ set $y_n=-x_n +2k_n\pi.1_A,$ with the convention $y_n=x_n-2k_n\pi.1_A$ when $\beta=0.$ Then $\chi(y_n)=n\beta=\chi(nx_1),$ and $cos(y_n)=cos(x_n)=c_n
=cos(nx_1).$ Since $\widehat A$ is a singleton, $y_n-n\beta.1_A$ and $nx_1-n\beta.1_A$ are quasinilpotent. When $\beta \neq  0$ and $\beta \neq \pi,$ it follows from item (i) of Lemma 2.2 that  $y_n-n\beta.1_A=nx_1-n\beta.1_A,$ and hence $y_n=nx_1.$

But $\Vert y_n-n\beta.1_A\Vert =\Vert y_n-\chi(y_n).1_A\Vert =\Vert \pm(x_n-\chi(\pm x_n).1_A)\Vert \le (1+\pi)M.$ Hence, given that $\Vert y_n-n\beta.1_A\Vert =n\Vert x_1-\beta.1_A\Vert$, we see that $x_1=\beta.1_A,$ $(c_n)_{n\ge 1}$ is a scalar sequence, and $c_n=cos(n\beta).1_A$ for $n \ge 1.$

If $\beta=0,$ then $\chi(x_n)\in [0,\pi]\cap2\pi\Z=\{0\},$ $y_n=x_n$ and $nx_1$ are quasinilpotent, and it follows from item (ii) of Lemma 2.2 that $x_n^2=n^2x_1^2.$ Since the sequence $(x_n^2)_{n\ge 1}$ is bounded, we have $x_1=0,$ $c_1=1_A$ and so $c_n=1_A$ for every $n \ge 1.$ If $\beta =\pi,$ set $c'_n=c_{2n}$ for $n\ge 1,$ and set $\beta'=arccos(\chi(c'_1)).$ Then $(c'_n)_{n\ge 1}$ is a cosine sequence. Since $\beta'=0,$ we have $2c_1^2-1_A=c_2=c'_1=1_A,$ and $(c_1-1_A)(c_1+1_A)=0.$ Since $scpec(c_1)=\{-1\},$ $c_1-1_A$ is invertible, $c_1=-1_A,$ and $c_n=(-1)^n.1_A=cos(n\pi).1_A$ for $n \ge 1.$ $\square$

Recall that a commutative unital Banach algebra $A$ is said to be local if it possesses a unique maximal ideal. We obtain the following corollary.

\begin{corollary} Let $G$ be an abelian group, and let $C= (C(g))_{g\in G}$ be a bounded  cosine family in a commutative unital local Banach algebra. Then $C$ is scalar, and so there exists a bounded cosine family $(c(g))_{g\in G}$ with values in $[-1,1]$ such that $C(g)=c(g).1_A$ for $g\in G.$

\end{corollary}

\section{When the distance to a bounded cosine function is strictly less than 2}

 A standard result shows that every bounded complex cosine function $c$ takes values in $[-1,1]$. The following observation, which is the cosine counterpart of a standard result for discontinuous one-parameter unimodular groups, see \cite{hp}, section 4.17, is certainly well-known.

\begin{proposition} Let $(c(t))_{t\in \R}$ be a discontinuous bounded complex cosine function. Then for every $\alpha \in [-1,1]$ there exists a sequence $(t_n)_{n\ge 1}$ of positive real numbers such that $lim_{n\to +\infty}t_n=0$ and $lim_{n\to+\infty}c(t_n)=\alpha.$
\end{proposition}
Proof: The identity

$$(1-c(s-t))(1-c(s+t))=(c(s)-c(t))^2$$

shows as is well-known  that this bounded  cosine function with values in $[-1,1]$  is discontinous at $0.$ Denote by $G$ the set of all real numbers $x$ for which there exists a sequence $(t_n)_{n\ge 1}$ of positive reals such that $lim_{n\to +\infty}t_n=0$ and $lim_{n\to +\infty}c(t_n)=cos(x).$ Then $nG +2\pi\Z\subset G$ for every $n\in Z,$ $G$ is closed, and $G\neq 2\pi \Z.$
Let $x \in G\cap (0,\pi].$ If $x\over \pi$ is irrational, then the sequence $(e^{inx})_{n\ge 1}$ is dense in $\T,$ and so $G=\R.$ If $x\over \pi$ is rational let $u$ be the smallest positive integer such that $e^{iux}=1.$ Then $e^{2i\pi\over u}=e^{ipx}$ for some $p\ge 1,$ and so ${2\pi\over u}\in G.$ Let $(t_n)_{n\ge 1}$ be a sequence of positive reals converging to 0 such that $lim_{n\to +\infty}c(t_n)=cos\left ({2\pi\over u}\right ),$ let $q\ge 2,$ and let $\alpha$ be a limit point of the sequence $c\left ({t_n\over u^{q-1}}\right )_{n\ge 1}.$ There exists $y \in \R$ such that $cos(y)=\alpha,$ so that $y \in G,$ and such that $y= {2\pi\over u^q} +{2k\pi\over u^{q-1}}=(1+ku){2\pi\over u^q}.$ Then $gcd(1+ku,u^q)=1,$ and there exist a positive integer $r$ such that $ry - {2\pi\over u^q} \in 2\pi\Z,$ so that  ${2\pi\over u^q}\in G.$ This implies that $G=\R.$ $\square$

\begin{corollary} Let $a\in \R,$ and let $(c(t))_{t\in \R}$ be a discontinuous bounded scalar cosine function. Then $sup_{t\in \R}\vert cos(at)-c(t)\vert=lim sup_{t\to 0}\vert cos(at)-c(t)\vert
=2.$
\end{corollary}

\begin{corollary} Let $X$ be a Banach space, let $(c(t))_{t\in R}$ be a scalar cosine function, and let $(C(t))_{t\in \R}$ be a bounded strongly continuous cosine family of bounded operators
on $X$ such that $sup_{t\in \R}\Vert C(t)-c(t)I_X\Vert <2.$ Then $c(t)$ is continuous, and so there exists $a\in \R$ such that $c(t)=cos(at)$ for $t\in \R.$
\end{corollary}

Proof: Let $x \in X$ such that $\Vert x \Vert=1.$ If $c(t)$ were discontinuous, there would exist a sequence $(t_n)_{n\ge 1}$ of positive real numbers such that $$lim_{n\to +\infty}t_n=lim_{n\to+\infty}c(t_n)+1=0,$$ which gives

$$sup_{t\in \R}\Vert C(t)-c(t)I_X\Vert \ge lim_{n\to +\infty}\Vert C(t_n)x -c(t_n)x\Vert=lim_{n\to +\infty}(1-c(t_n))\Vert x \Vert=2.$$

$\square$

The following observation is an easy consequence of Kronecker's theorem on independent finite subsets of the unit circle.

\begin{lemma} Let $a,b$ be two real numbers. If $pa-qb\neq 0$ for $(p,q) \in \Z^2\setminus \{0,0\},$ then $sup_{t\in \R}\vert cos(ta)-cos(tb)\vert =2.$ 

\end{lemma}

Since  $sup_{t \in \R}\vert 1-cos(bt)\vert=2$ for every $b\neq 0,$ we can restrict attention to the case where $a\neq 0$ and $b\neq 0.$ Denote by $\T$ the unit circle. If $pa-qb \notin 2\pi \Z$ for $(p, q)\in \Z^2\setminus \{(0,0)\},$ then the set $\{(e^{ia},e^{ib})\}$ is independent, and it follows from Kronecker's theorem, see for example \cite{ks}, p. 21 that the set $\left \{(e^{ina},e^{inb})\right \}_{n\in \Z}$ is dense in $\T^2,$  which implies that $sup_{n\ge 1}\left \vert cos(na)-cos(nb)\right \vert =2.$ 

Now assume that $pa -qb =2k\pi$ for some $(p,q)\in  \Z^2\setminus \{(0,0)\}$ and some $k\in \Z.$

If $p=0,$ then we have $b={2k\pi \over q},$ with $k\ge 1,$ $q\ge 1.$ If $a\over \pi$ is irrational, then $\{e^{inqa}\}_{n\ge 1}$ is dense in $\T,$ and we have

$$sup_{n\ge 1}\left | cos(naq) -cos(nqb)\right | = sup_{n\ge 1}\left\vert cos(naq)-1\right \vert=2.$$

Otherwise $a \over \pi$ is rational, and so is ${a\over b}.$

Now assume that $p\neq 0, q\neq 0,$ and $k\neq 0.$ If $a\over \pi$ is rational, then $b\over \pi$ is rational, and ${a\over b}$ is rational. Otherwise we have,  since the set $\{e^{ian}\}_{n\ge 1}$ is dense in the unit circle,

$$sup_{t\in \R}\left | cos(bt) -cos(at)\right  | \ge sup_{n\ge 1} \left | cos\left ( b(2kn+1){q\over 2k} \right ) -cos \left( a(2kn+1){q\over 2k} \right)\right |$$

$$=\left | cos\left (a (2kn+1){p\over 2k}  -(2kn+1)\pi \right )- cos \left( a(2kn+1){q\over 2k} \right)\right |$$ $$=\left | cos\left ( \left (an+{a\over 2k}\right ){p}  \right )+ cos \left( \left (an+{a\over 2k}\right ){q} \right)\right | =2.$$

Hence $k=0$ if $sup_{t\in \R}\vert cos(ta)-cos(tb)\vert <2.$ $\square$

\begin{lemma} Let $a \ge 0,$ and set \  $\Omega(a,m):=\left \{b \ge 0 \ | \ sup_{t\in \R} \left | cos(ta)-cos(tb) \right \vert \le m\right\}$ for $m\ge 0.$

(i) If $m<2,$ then $\Omega(a,m)$ is finite, and every $b\in \Omega(a,m)$ has the form $b={pa\over q},$ where $p,q$ are odd, $gcd(p,q)=1$, $1\le p \le {\pi\over arccos(m-1)}, 1\le q \le {\pi\over arccos(m-1)}.$

(ii)  If $m<{8\over 3\sqrt 3},$ then $\Omega(a,m)=\{a\}.$

\end{lemma}

Proof: Assume that $m<2.$ Since $\Omega(0,m)=\{0\},$ we can assume that $a > 0.$ Let $b\in \Omega(a,m).$ Then $b \neq 0.$ We may restrict attention to the case where $a=1$, and  there exists positive integers $p$ and $q$ such that $b={p\over q},$ with $gcd(p,q)=1.$ If $p$ or $q$ were even, we would have $sup_{t\in \R}\vert cos(t)-cos(bt)\vert \ge sup_{n\ge 1}\vert cos(np\pi)-cos(nq\pi)\vert =2,$ and so $p$ and $q$ are odd.

 It follows from Bezout's theorem that there exists $(u,v) \in \Z$ such that $2up-(2v+1)q=1.$ We have

$$sup_{t\in \R}\vert cos(t)-cos(bt)\vert \ge cos(2u\pi)-cos\left ({2up\pi\over q}\right )=1- cos\left ((2v+1)\pi +{\pi\over q} \right )$$ $$=1+cos\left ({\pi\over q}\right ).$$

Hence $q \le {\pi\over arccos(m-1)}.$ Since $sup_{t\in \R} \vert cos(t)-cos(bt)\vert =sup_{t\in \R} \vert cos(t)-cos({t\over b})\vert,$ the same argument shows that $p\le {\pi\over arccos(m-1)}.$ This proves (i)

\smallskip
(ii) If $m< 1+cos\left ({\pi\over 5}\right )\approx 1,8090,$ then  $arccos(m-1)>{\pi\over 5}$ and so every  $b\in \Omega(a,m)$ can be written under the form $b={ap\over q},$ where $1\le p <5$ and $1\le q <5,$ with $p$ and $q$ odd, $gcd(p,q)=1$ and so $\Omega(a,m)\subset \left \{{a\over 3}, a, 3a \right \}.$

 An elementary computation shows that $\vert cos(t) -cos(3t)\vert $ attains its maximum when  $cos(t) =\pm {1\over \sqrt 3}$, which gives
 $$max_{t \in \R}=\vert cos(t) -cos\left ({t/ 3}\right )\vert =max_{t \in \R}\vert cos(t) -cos(3t)\vert = {8\over 3\sqrt 3}\approx 1.5396.$$ Hence 
  $\Omega(a,m)=\left \{{a\over 3}, a, 3a \right \}$ if ${8\over 3\sqrt 3}\le m < 1+cos\left ({\pi\over 5}\right ),$ and   $\Omega(a,m)=\{ a\}$ if $m<{8\over 3\sqrt 3}.$

$\square$

We obtain the following theorem.

\begin{theorem} Let $(C(t))_{t\in \R}$ be a cosine function in a unital Banach algebra, and let $A_1$ be the closed subalgebra of $A$ generated by $(C(t))_{t\in \R}.$

(i) If $m:=sup_{t\in \R}\left \Vert  C(t)-cos(at).1_A\right\Vert <2$ for some $a\in \R,$ then $A_1$ is isomorphic to $\C^k$ for some $k\le card\left ( \Omega(a,m)\right ),$ and there exists a family $p_1,\dots,p_k$ of pairwise orthogonal  idempotents of $A_1$ and a family $b_1,\dots, b_k$  of distinct  elements of $\Omega(a,m)$  such that $C(t)=\sum \limits_{j=1}^kcos(b_jt)p_j$ for $t \in \R.$ 

(ii)  If $sup_{t\in \R}\left \Vert  C(t)-cos(at).1\right\Vert <{8\over 3\sqrt 3}$ for some $a\in \R,$ then $C(t)=cos(at).1_A$ for $t\in \R.$
\end{theorem}

Proof: (i) Let $\chi \in \widehat{A_1}.$ Then it follows from Corollary 3.2 that there exists $b_{\chi} \in \Omega(a,m)$ such that $\chi(C(t))=cos\left (tb_{\chi}\right )$ for $t\in \R.$ 

Since $b^2_{\chi}=2 lim_{t\to 0}{1-\chi(C(t))\over 2},$ the map $\chi \to b_{\chi}$ is one-to-one, and it follows from Lemma 3.5 that $\widehat A_1$ is finite.

 Let $\chi_1, \dots, \chi_k$ be the elements of $\widehat{A_1}.$ It follows from the standard one-variable holomorphic functional calculus, see for example \cite{d}, that there exists for every $j\le k$ an idempotent $p_j$ of $A_1$ such that $\chi_j(p_j)=1$ and $\chi_i(p_j)=0$ for $i\neq j.$ Hence $p_jp_i=0$ for $j \neq i,$ and $1_A=1_{A_1}=\sum \limits _{j=1}^kp_j.$ 

Let $x\in A_1.$ Then $(p_jC(t))_{t\in \R}$ is a cosine function in the commutative unital Banach algebra $p_jA_1,$ $spec_{p_jA_1}(p_jC(1))=\{\chi_j(C(1))\},$ and $(cos(at)p_j)_{t\in \R}$ is a scalar cosine function in $p_jA_1.$ Since $sup_{t\in \R}\left \Vert cos(at) p_j-p_jC(t) \right \Vert \le 2 \Vert p_j\Vert,$ the cosine function $(p_jC(t))_{t\in \R}$ is bounded, and it follows from theorem 2.3  that $(p_jC(nt))_{n\in \Z}$ is a scalar cosine sequence for every $t\in \R.$ So $(p_jC(t))_{t\in \R}$ is a scalar cosine function, and $p_jC(t)=\chi_j(C(t)p_j=cos(b_jt)p_j,$ where $b_j=b_{\chi_j}\in \Omega(a,m).$

We obtain

$$C(t)=\sum \limits_{j=1}^kC(t)p_j=\sum \limits_{j=1}^kcos(tb_j)p_j \ (t\in \R).$$

Since the algebras $p_jA_1$ are one-dimensional, $A_1$ is isomorphic to $\C^k.$

(ii) If $m<{8\over 3\sqrt 3},$ then $\Omega(a,m)=\{a\},$ $k=1,$ $p_1=1_A$ and $C(t)=cos(at).1_A.$ $\square$

Since $\Omega(0,m)=\{0\}$ for every $m<2,$ we obtain the following result, which was obtained recently by Schwenninger and Zwart in \cite{zs1} for strongly continuous operator valued cosine functions. A very short argument to prove a weaker result with the constant $3/2$ instead of $2$ is given by Arendt in \cite{a}.

\begin{corollary} Let $(C(t))_{t\in \R}$ be a cosine function in a unital Banach algebra $A.$ If $sup_{t\in \R}\Vert C(t)-1_A\Vert <2,$ then $C(t)=1_A$ for $t\in \R.$
\end{corollary}

 Let $C=(C(t))_{t\in \R}$ be a bounded strongly continuous cosine family of bounded operators
on a Banach space $X.$ Bobrowski and Chojnacki observed in \cite{bc}, lemma 3 that if $sup_{t\in \R}\Vert C(t)-c(t)I_X\Vert <1$ then the generator of $C$ is bounded, assuming that the scalar cosine function $c$ is continuous. The following corollary, which is an immediate consequence of collorary 3.3 and theorem 3.6, shows that the continuity condition on $c(t)$ is redundant, and that the generator of $C(t)$ is bounded whenever $sup_{t\in \R}\Vert C(t)-c(t)I_X\Vert <2.$
\begin{corollary}  Let $X$ be a Banach space, let $(c(t))_{t\in R}$ be a bounded scalar cosine function, and let $C=(C(t))_{t\in \R}$ be a bounded strongly continuous cosine family of bounded operators
on $X$ such that $sup_{t\in \R}\Vert C(t)-c(t)I_X\Vert <2.$ Then there exists $a\in \R$ such that $c(t)=cos(at)$ for $t\in \R,$ and the conclusions of theorem 3.6 hold. In particular the generator of the cosine function $C$ is bounded.
\end{corollary}

{\it Jean Esterle

IMB, UMR 5251

Universit\'e de Bordeaux

351, cours de la Lib\'eration

33405 -Talence (France)

esterle@math.u-bordeaux1.fr}

\end{document}